\documentclass[11pt,a4paper]{article}
\usepackage[utf8]{inputenc}

\usepackage{authblk}
\usepackage{geometry}
\usepackage{array}
\usepackage{booktabs}
\usepackage{multirow}
\usepackage{lipsum}
\usepackage{color}
\usepackage{amsmath}
\usepackage{amsthm}
\theoremstyle{remark}

\usepackage{amsfonts}
\usepackage{amssymb}
\usepackage{bm}
\usepackage{caption}
\usepackage{overpic}  
\usepackage{graphicx,subfig}
\usepackage[colorlinks,linkcolor=blue,anchorcolor=red,citecolor=blue]{hyperref}
\usepackage{rotating}
\usepackage{tikz}
\usetikzlibrary{positioning, arrows.meta}
\usepackage{algorithm, algorithmic}
\usepackage{cite}

\usepackage{setspace}
\usepackage{multirow}
\numberwithin{equation}{section}  
\newtheorem{thm}{Theorem}[section]

\newtheorem{rmk}[thm]{Remark}

\def\d{\displaystyle}

\title{Optimal control of stochastic cylinder flow using data-driven compressive sensing method}
\author[a]{Liuhong Chen}
\author[a]{Ju Ming\thanks{Corresponding author: jming@hust.edu.cn}}
\author[b]{Max D. Gunzburger}
\affil[a]{School of Mathematics and Statistics, Huazhong University of Science and Technology, Wuhan 430074, China}
\affil[b]{Department of Scientific Computing, Florida State University, Tallahassee, USA}

\date{} 

\begin{document}
      \newtheorem{theorem}{Theorem}[section]
      \newtheorem{assumption}[theorem]{Assumption}
      \newtheorem{corollary}[theorem]{Corollary}
      \newtheorem{proposition}[theorem]{Proposition}
      \newtheorem{lemma}[theorem]{Lemma}
      \newtheorem{definition}[theorem]{Definition}
      \newtheorem{algo}[theorem]{Algorithm}
      \newtheorem{remark}[theorem]{Remark}
      
      \maketitle
\noindent \rule[4pt]{17cm}{0.05em}  

\noindent   
\textbf{Abstract:} 
A stochastic optimal control problem for incompressible Newtonian channel flow past a circular cylinder is used as a prototype optimal control problem for the stochastic Navier-Stokes equations. The inlet flow and the rotation speed of the cylinder are allowed to have stochastic perturbations. The control acts on the cylinder via adjustment of the rotation speed. Possible objectives of the control include, among others, tracking a desired (given) velocity field or minimizing the kinetic energy, enstrophy, or the drag of the flow over a given body. Owing to the high computational requirements, the direct application of the classical Monte Carlo methods for our problem is limited. To overcome the difficulty, we use a multi-fidelity data-driven compressive sensing based polynomial chaos expansions (MDCS-PCE). An effective gradient-based optimization for the discrete optimality systems resulted from the MDCS-PCE discretization is developed. The strategy can be applied broadly to many stochastic flow control problems. Numerical tests are performed to validate our methodology. \\


\noindent \rule[4pt]{17cm}{0.05em} 

\maketitle

\section{Introduction}\label{sec1}

Stochastic Navier-Stokes (SNS) equations have broad applications in diverse fields such as atmospheric science, aero- and oceanographic engineering, etc.; see e.g., \cite{Breit2015Stochastic, Constantin2014Unique, Frisch2010FULLY, Mason2010Large, menaldi2002stochastic, Shafiq2017Statistical, SILAEN2010Effect}. In recent years, they have been extensively used to simulate the random flow motion of incompressible Newtonian fluids. However, due to the high computational requirements, numerical simulations of control problems involving stochastic flows described by the SNS equations still remain big challenges from the computational point of view. The classical Monte Carlo (MC) method \cite{metropolis1949monte, Rubinstein2008Simulation} can be used to determine reliable statistical moments of the solution, thanks to the \emph{Central Limit Theorem} (CLT). However, to obtain acceptable accuracy that meets engineering specifications in technological applications, one has to balance the statistical errors governed by CLT and the numerical errors for each sample \cite{Giles2008Multilevel}, thus often resulting in a prohibitively high computational complexity. The latter is largely due to the requirement of the large size of the sample solutions to minimize the statistical error and fine spatial mesh to accurately determine each sample solution. Additionally, if a gradient-based optimization method is applied, substantial additional computational effort is required to deal with the solution of the associated adjoint system. To avoid these dilemmas and motivated by the successes of polynomial chaos expansions (PCEs) as a non-statistical approach for stochastic partial differential equations (SPDEs) (see e.g., \cite{Dolgov2015Polynomial, ghanem2003stochastic, hou2006wiener,  sakamoto2002polynomial, Sepahvand2010UNCERTAINTY, xiu2003modeling, Karagiannis2014, Contreras2018, Yang2019}), in the present work, we study a new multi-fidelity data-driven compressive sensing based PCEs (MDCS-PCE) for stochastic flow control problems, using a boundary control problem of stochastic cylinder flow as an illustration. We refer to \cite{borzi2009multigrid, gunzburger2011optimal, gunzburger2003perspectives, gunzburger2018, Ali2017, Kouri2013, Tiesler2012} and references therein for the multi-fidelity modeling and optimal control problems constrained by PDEs with random inputs.

PCEs \cite{Dolgov2015Polynomial, ghanem2003stochastic, hou2006wiener, sakamoto2002polynomial, Sepahvand2010UNCERTAINTY, xiu2003modeling}, originated from the idea of \emph{homogeneous chaos} introduced by Wiener \cite{Furth1960Non, Wiener1},  was applied by Ghanem and Spanos \cite{sakamoto2002polynomial, Ghanem1991Stochastic} as a means to resolve a variety of mechanics problems involving random inputs. As a non-statistical method for SPDEs, it is very useful when the number of uncertain parameters is moderate. In this framework, a given random variable is represented as a series in terms of a suitable orthogonal basis, such as the Hermite polynomials related to the Gaussian random variables \cite{Ernst2011On, xiu2002wiener}. After establishing the Galerkin approximation via PCEs \cite{babuvska2007stochastic, Teckentrup2014A}, the primary task is to determine the PCE coefficients by solving the resulted intrusive systems, i.e., the coupled system of deterministic PDEs with the same type of nonlinearity as the original SPDEs. Obviously, the computational complexity is proportional to the number of polynomial chaos basis in the intrusive system with the size $\mathcal{O}{(N\times n)}$, where $N$ is the number of polynomial chaos basis and $n$ is the nodes, which would lead to the \emph{curse of dimensionality} when dealing with the uncertainty problems in high-dimensional probability space. Fortunately, due to the fact that many high-dimensional problems have intrinsic sparsity property, i.e., the number of the nonzero or large coefficients of the PCEs is generically much smaller than the cardinality of the basis, we are able to incorporate the compressive sensing (CS) method (\cite{Arjoune2018Compressive, baraniuk2007compressive, candes2006stable, candes2008introduction, candes2006compressive, Zhang2017A, Doostan2011, BRUCKSTEIN2009}) emerged from the field of sparse signal recovery. There have been lots of works using CS to solve SPDEs, some of the relevant analyses can be found in \cite{Bouchot2015Compressed, doostan2011non, jakeman2015enhancing, mathelin2012compressed, peng2014weighted, Rauhut2015Compressive}.

  The MDCS-PCE method studied in this work utilizes the development of similar ideas in our earlier work \cite{Liang}. The essential idea of MDCS-PCE is to construct problem-dependent bi-orthogonal bases in the expansion of the sample solutions, in which the coefficients become more sparse. A multi-fidelity approach is used to do data generation on a coarse mesh upon which the statistical information based on the CS method can be used for computation on a fine mesh. Thus, to achieve the same accuracy as the classical CS method, MDCS-PCE needs far less measurements, which offers the potential to make the MDCS-PCE method very effective in practice. Under the hypothesis that the solutions of the stochastic constraint functions are reasonably compressible, MDCS-PCE requires much less number of calls of the deterministic solver than other sampling methods such as MC or sparse grid collocation. Since the cost of solving the deterministic problem is always the main cost of all, the MDCS-PCE method greatly reduces the computational cost for solving the stochastic optimal control problem, and is more efficient for high-dimensional problems. The numerical results in this work provide further confirmation.

  The rest of the paper is organized as follows. In Section \ref{sec2}, we present the model problem used for our study and the noise discretization schemes we employ. We also give the weak formulation and the specific discrete form of the stochastic cylinder flow problem. Then in section \ref{sec3}, we briefly describe the general framework of the MDCS-PCE method for solving the stochastic cylinder flow problem together with some numerical results. In Section \ref{sec4}, we formulate the stochastic boundary control problems for the stochastic cylinder flow, then combined with gradient descent iteration, an optimization algorithm based on the MDCS-PCE method is proposed, along with some computational experiments to validate our optimization algorithm. Finally, in Section \ref{sec5}, some concluding remarks are given. 

\section{Stochastic cylinder flow}\label{sec2}

\begin{figure}[h!]
   \begin{center}
   \includegraphics[width=5in]{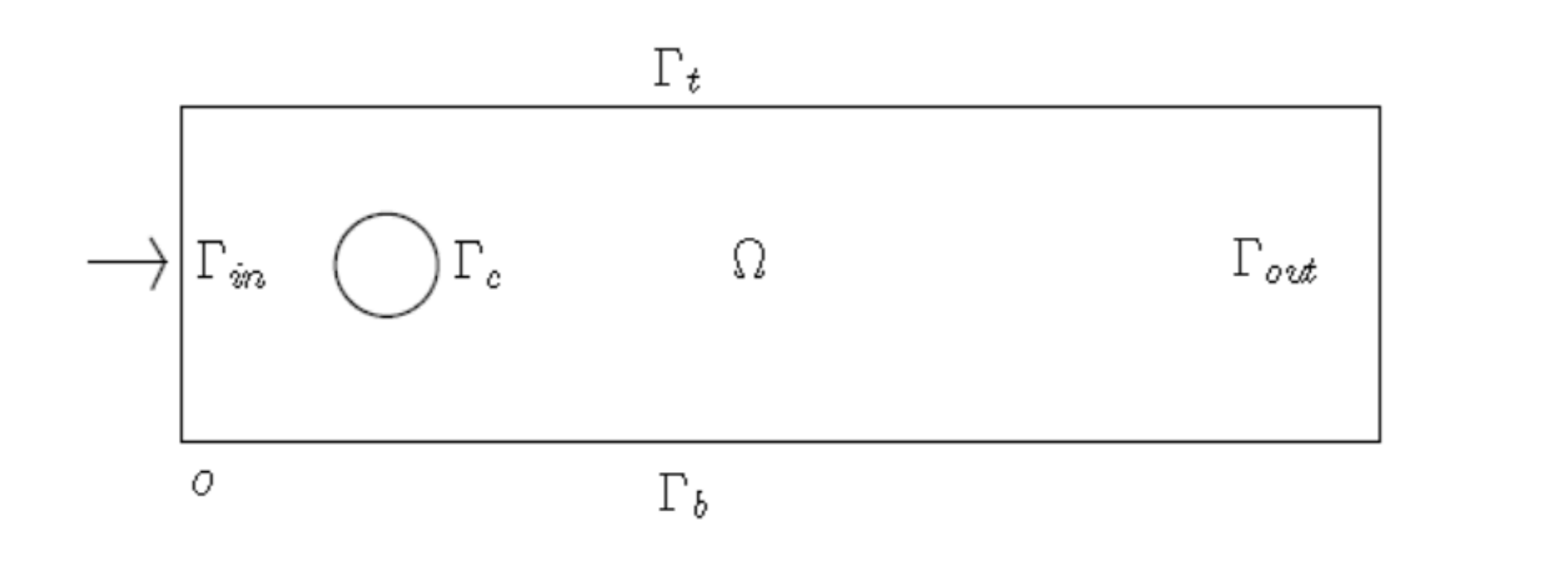}
   \caption{\textnormal{\textnormal{Channel flow past a circular cylinder.}}}
   \label{cylinder_domain}
   \end{center}
\end{figure}

In this section, we briefly describe the stochastic Navier-Stokes (SNS) equations for modeling incompressible Newtonian flow and the associated numerical schemes for generating PCE solutions. To give the discussion context, we consider the concrete problem of the two-dimensional stochastic incompressible flow past a circular cylinder of diameter $d$ over the time interval $[0,T]$. We denote by $\Omega$ the channel domain, whose boundary $\Gamma$ consists of four parts as depicted in Figure \ref{cylinder_domain}. Furthermore, we choose the origin $o$ to be at the bottom-left corner of the channel and $(x_{0},y_{0})$ to denote the center of the circular cylinder.

For $t\in(0,T]$ and $\textbf{x}\in\Omega$, the flow is governed by the non-dimensionalized Navier-Stokes equations
\begin{gather}
   \textbf{u}_{t}-\frac{1}{Re}\Delta \textbf{u}+ (\textbf{u}\cdot \nabla)\textbf{u}+\nabla p=\textbf{0} \qquad \text{ in } \Omega,\\
   \nabla\cdot \textbf{u}=0  \qquad \text{ in } \Omega,
   \label{NS}
\end{gather}
where $\textbf{u}=(u,v)^{T}$ denotes the velocity field, $p$  the pressure, $Re=\rho UH/\mu$ the usual Reynolds number, $\rho$ the density of the fluid, $U$ is some measure of the average speed of the inflow, $H$ the channel height, and $\mu$ the kinematic viscosity. The problem specification is completed by the imposition of an initial condition on the velocity and the boundary conditions, for $t\in(0,T]$,
\begin{gather}
   \textbf{u}|_{\Gamma_{in}}=(u_{in},0)^{\top},\label{bond1}\\
   \textbf{u}|_{\Gamma_{t}\cup\Gamma_{b}}=(0,0)^{\top},\label{bond2}\\
   \textbf{u}|_{\Gamma_{c}}=\big(-(y-y_{0})\phi, (x-x_{0})\phi\big)^{\top},\label{gammac}\\
   \left(p\textbf{n}-\frac{1}{Re}\frac{\partial \textbf{u}}{\partial \textbf{n}}\right)\big|_{\Gamma_{out}}=(0,0)^{\top},
   \label{donothingcond}
\end{gather}
where $\phi$ denotes the angular velocity of the cylinder. Due to the fact that the velocity field $\textbf{u}$ on $\Gamma_{out}$ is in general not known, here, the outflow condition \eqref{donothingcond} is "{\it artificially}'' imposed; it has been found to be an acceptable approximation, provided that the outflow boundary is placed at sufficient distance downstream of the cylinder; see e.g., \cite{Sani2010R} for further discussion. Other conditions, e.g., $p=0$, may also be imposed along $\Gamma_{out}$. However, we do note that the condition \eqref{donothingcond} may fail for flows with high Reynolds numbers; see, e.g., \cite{Fehn2018Robust,turek1999efficient}.

In this paper, the stochastic inputs act on the inflow condition $u_{in}$ and the angular velocity $\phi$. Specifically, we let
\begin{eqnarray}
   u_{in}(0, y, t; \bm{\xi})=\overline u_{in}(0,y,t) +\sigma_{1}\dot \omega(t;\bm{\xi})  \quad \text{ and }\quad \phi(t;\bm{\xi})=\overline{\phi}(t)+\sigma_{2}\dot \omega(t;\bm{\xi}),
   \label{sto_cond}
\end{eqnarray}
where $\omega(t)$ represents \emph{Brownian} noise and the constants $\sigma_{i}\ge 0$, $i=1,2$, and functions $\overline u_{in}$ and $\overline{\phi}(t)$ are deterministic.

\subsection{Representation of white noise}\label{sec2.1}

It is well known that the \emph{Brownian} motion ${\omega}(s)$ can be approximated via Fourier expansion. Particularly, letting $\{\varphi_{i}\}_{i=1}^{\infty}$ denote an orthonormal basis of $L^2([0,t])$, we have
\begin{eqnarray}
   {\omega}(s)=\sum_{i=1}^{\infty} \xi_{i}\int_{0}^{s}\varphi_{i}(\tau) d\tau ,
   \label{rep_Brownian}
\end{eqnarray}
where $\xi_i\sim\mathcal N(0,1)$ are independent and identically distributed ($i.i.d.$) random variables. In fact, one can easily see that the It\^o integral $\int_{0}^{t}\varphi_{i}d\omega(s) \stackrel{\text{i.i.d.}}{\sim} \mathcal N(0,1)$. Thus, if we let
$$
\xi_{i}=\int_{0}^{t}\varphi_{i}d\omega(s),
$$
then
\begin{equation}
   \begin{aligned}
   \omega(s) = \int_{0}^{t} \chi_{[0,s]}(\tau)d\omega(\tau) & = \sum_{i=1}^\infty\int_{0}^{t}(\chi_{[0,s]},\varphi_{i})\varphi_{i}(\tau)d\omega(\tau)\\
   &\qquad\quad = \sum_{i=1}^{\infty} \xi_{i}\int_{0}^{s}\varphi_{i}(\tau)d\tau,
   \end{aligned}
   \label{noise_exp}
\end{equation}
where $\chi_{[0, s]}(\tau)$ is the characteristic function of the interval $[0, s]$. It is well known that the expansion \eqref{noise_exp} uniformly converges to $\omega(s)$ in the mean square sense, i.e., $\forall
s\le t$,
$$
\lim_{N\to\infty}\mathbb{E}\Big[\omega(s)-\sum_{i=1}^{N} \xi_{i}\int_{0}^{s}\varphi_{i}(\tau)d\tau\Big]^2=0.
$$
In \cite{du2002numerical}, comparisons between noises generated by different bases $\{\varphi_{i}\}_{i=1}^{\infty}$ are provided.

\subsection{Polynomial chaos expansion}\label{sec2.2}

To start with, define the following set of multi-indices with the finite number of non-zero components:
$$
\mathcal{I}:=\Big\{\bm{\alpha}=(\alpha_{1}, \alpha_{2}, \cdots)\,:\, \alpha_{i}\in \mathbb N^{+}_{0}, \,|\bm{\alpha}|=\sum_{i=1}^{\infty} |\alpha_{i}|<\infty\Big\}. 
\label{multi-indices}
$$
The $n$-th order normalized one-dimensional Hermite polynomials are defined as
$$
h_{i}(x)=(-1)^{i}e^{x^2/2}\frac{d^{i}}{d x^{i}}e^{-x^{2}/2},\qquad i=0, 1, 2, \ldots.
$$
It is well known that the set  $\{h_{i}(x)\}_{i=1}^{\infty}$ is a complete orthonormal basis of $L^{2}( \mathbb{R})$ with respect to (w.r.t.) Gaussian weighting function $\rho$ (Gaussian measure), i.e.,
$$
\int_{\mathbb R}\rho(x)h_{i}(x)h_{j}(x)dx= \mathbb{E}[h_{i}h_{j}] =\delta_{i j},
$$
where
$
\rho(x)=\frac{1}{\sqrt{2\pi}}e^{-x^2/2}
$
denotes the standard Gaussian distribution function and $\delta_{i j}$ the Kronecker delta function.  The one-dimensional  Gaussian space is then defined as
$$
L^2_\rho(\mathbb{R}):=\Big\{f(x): \int_\mathbb{R}\rho(x)|f(x)|^2dx<\infty\Big\}.
$$
Clearly, tensor products of the elements of $\{h_{i}\}_{i=1}^{\infty}$ construct a complete basis of the corresponding multi-dimensional Gaussian probability space. Specifically, for a multi-index $\bm{\alpha}=(\alpha_{1},\alpha_{2}, \cdots)\in \mathcal{I}$, a polynomial chaos basis functions of order $|\bm{\alpha}|$ is defined as
\begin{gather}
   \mathcal{H}_{\bm{\alpha}}(\bm{\xi}):=\prod_{i=1}^{\infty} h_{\alpha_{i}}(\xi_{i}),
   \label{Wick}
\end{gather}
where $|\bm{\alpha}|=\sum_{i=1}^{\infty} |\alpha_{i}|$, $\bm{\xi}=(\xi_{1}, \xi_{2},\cdots)$. We have the following fundamental theorem for PCE methods.

\begin{theorem}\label{CMtheorem}{\rm[Cameron--Martin theorem \cite{cameron1947orthogonal}]}
   For any given point $(\bm x, s)\in \Omega\times [0,t)$, if  $u(\bm x,s)$ is a stochastic functional of Brownian motion $\omega$ with $\mathbb{E}[u(\bm x,s)]^{2} < \infty$, then $u(\bm x,s)$ can be represented by the expansion
   \begin{gather}
      u(\bm x,s)=\sum_{\bm{\alpha}\in \mathcal{I}}u_{\bm{\alpha}}(\bm x,s)\mathcal{H}_{\bm{\alpha}}(\bm{\xi}),
      \label{PCE}
   \end{gather}
   where $\{\mathcal H_{\bm{\alpha}}\}_{\bm{\alpha}\in\mathcal I}$ is the set of polynomials defined by \eqref{Wick} and $\bm{\xi}=(\xi_{1},\xi_{2},\cdots)$ with i.i.d. random variables $\xi_i\sim\mathcal N(0,1)$; $u_{\bm{\alpha}}$, $\bm{\alpha}\in\mathcal I$, are referred to as the \emph{chaos coefficients.} Furthermore, the first two statistical moments of $u(\bm x, s)$ are given by
   $$
   \mathbb{E}[u(\bm x, s)]=u_{0}(\bm x, s ) \quad\text{ and }\quad \mathbb{E}[u(\bm x,	s)^2]=\sum_{\bm{\alpha}\in \mathcal{I}}|u_{\bm{\alpha}}(\bm x, s)|^2.
   $$
\end{theorem}

\begin{rmk}\label{rmk}
   From the definition of the  polynomials \eqref{Wick} and the representation of Brownian motion \eqref{rep_Brownian}, we have the following relationship between $\omega$ and the multi-dimensional Gaussian variable $\xi$:
   $$
   \omega(s) = \sum_{\stackrel{\bm{\alpha} \in \mathcal {I}}{|\bm{\alpha}| =1, \alpha_{j}=1}} \int_{0}^{s}\varphi_{j}(\tau)d\tau \,\mathcal H_{\bm{\alpha}}(\bm{\xi}).
   $$
\end{rmk}

\subsection{Formulation of the cylinder flow problem}\label{sec2.3}

To give the weak formulation of the system \eqref{NS}, we first define the stochastic Sobolev spaces 
\begin{equation}
   \mathcal{L}^{2}(D;0,T;H^{1}(\Omega))=L^{2}(0,T;H^{1}(\Omega)) \otimes L^{2}(D) ,     
   \label{sto_space}
\end{equation}
which is a tensor product space, here $D$ is the stochastic variable space, and $L^2(0,T;H^{1}(\Omega))$ is the space of strongly measurable maps $\mathbf{v}: [0,T] \rightarrow H^{1}(\Omega)$ such that
\begin{equation}
   \|\mathbf{v}\|^{2}_{L^{2}(0,T;H^{1}(\Omega))}=\int_{0}^{T} \|\mathbf{v}\|^{2}_{\Omega} dt ,
   \label{time_space}
\end{equation} 
so the norm of the tensor Sobolev space \eqref{sto_space} can be defined as 
\begin{equation}
   \| \mathbf{v}\|^{2}_{\mathcal{L}^{2}(D;0,T;H^{1}(\Omega))}=\int_{D} \|\mathbf{v}\|^{2}_{L^{2}(0,T;H^{1}(\Omega))} dP = \mathbb{E}\left [ \|\mathbf{v}\|^{2}_{L^{2}(0,T;H^{1}(\Omega))} \right ].
\end{equation} 
Then, a weak formulation of the NS equations \eqref{NS} is given as follows: given an initial velocity $\textbf{u}_0(\bm{x})\in L^2(\Omega),$  we seek $\textbf{u}\in \mathcal{L}^2(D;0,T;H^{1}(\Omega))$, such that \eqref{bond1}--\eqref{donothingcond} are satisfied and
\begin{gather}
   \mathbb{E}\left \{ \left ( \frac{\partial \textbf{u}}{\partial t},\textbf{v} \right) + a(\textbf{u},\textbf{v}) + c(\textbf{u},\textbf{u},\textbf{v}) + b(\textbf{v},p) \right \}=0 \nonumber \\
   \mathbb{E}[b(\textbf{u},q)]=0 \label{weak_form}\\
   \textbf{u}(0,\bm{x})=\textbf{u}_0(\bm{x})\nonumber
\end{gather}
for all test functions $ \textbf{v}\in V_{0} =: \{\textbf{w}\in \mathcal{L}^2(D;0,T;H^{1}(\Omega)): \textbf{w}|_{\Gamma\backslash\Gamma_{out}}=0\}$ and $q \in \mathcal{L}^2(D;0,T;L^{2}(\Omega))$; here, $(\cdot,\cdot)$ denotes the $L^{2}(\Omega)$ inner product. The details about the bilinear operators $a(\textbf{u},\textbf{v})$, $b(\textbf{v},p)$ and the trilinear operator $c(\textbf{u},\textbf{u},\textbf{v})$ can be found in \cite{gunzburger2003perspectives}.

Then by substituting \eqref{PCE} into \eqref{weak_form}, the associated weak formulation for the chaos coefficients is then given as follows: given an initial velocity $\textbf{u}_{0}(\bm{x})\in L^{2}(\Omega),$  we seek $\textbf{u}_{\bm{\alpha}}(\bm x, t)\in L^{2}(0,T;H^{1}(\Omega))$, $\forall \bm{\alpha} \in \mathcal{I}$, such that \eqref{bond2}--\eqref{donothingcond} are satisfied and
\begin{gather}
   \left ( \frac{\partial\textbf{u}_{\bm{\alpha}}}{\partial t},\textbf{v} \right )+a(\textbf{u}_{\bm{\alpha}},\textbf{v})+ \sum_{\bm{\beta}, \bm{\gamma} \in \mathcal{I}}C_{\bm{\beta}, \bm{\gamma},\bm{\alpha}} c(\textbf{u}_{\bm{\beta}},\textbf{u}_{\bm{\gamma}},\textbf{v}) + b(\textbf{v},p_{\bm{\alpha}}) = \textbf{0}\nonumber\\
   b(\textbf{u}_{\bm{\alpha}},q)=0, \label{Wick_weak1}\\
   \textbf{u}_{\bm{\alpha}}(0,\bm{x};\bm{\xi})=\textbf{u}_{0,\bm{\alpha}}(\bm{x};\bm{\xi})\nonumber
\end{gather}
for all test functions $ \textbf{v}\in V_{0}$ and $q \in L^{2}_{0}(\Omega)$, where
$$
C_{\bm{\beta},\bm{\gamma},\bm{\alpha}}=\frac{ \mathbb{E}[\mathcal{H}_{\bm{\beta}} \mathcal{H}_{\bm{\gamma}} \mathcal{H}_{\bm{\alpha}}] }{\mathbb{E}[\mathcal{H}_{\bm{\alpha}}\mathcal{H}_{\bm{\alpha}}]},\qquad \textbf{u}_{\alpha}=(u_{\bm{\alpha}}, v_{\bm{\alpha}})^{\top}.
$$
Note that
$$
\mathbb{E}[\mathcal{H}_{\bm{\alpha}}\mathcal{H}_{\bm{\beta}}]=\delta_{\bm{\alpha}\bm{\beta}}.
$$
Clearly,$C_{\bm{\beta},\bm{\gamma},\bm{\alpha}}=C_{\bm{\gamma},\bm{\beta},\bm{\alpha}}.$ 

The PCE solution is then determined by
\begin{eqnarray}
   \d\textbf{u}(\bm x, t;\bm{\xi})=\sum_{\bm{\alpha}\in\mathcal{I}}\textbf u_{\bm{\alpha}}(\bm x, t)\mathcal{H}_{\bm{\alpha}}(\bm{\xi}).
   \label{solution_expan}
\end{eqnarray}
We select finite element methods for spatial discretization and modified Newton linearization for time discretization \cite{gunzburger2011optimal}. Note that \eqref{Wick_weak1} is a coupled system of \emph{deterministic} PDEs for the chaos coefficients, solving this problem directly is complicated in programming, so we consider using our data-driven compressive sensing method instead to derive the PCE coefficients above.

\section{Data-Driven compressive sensing approach for stochastic cylinder flow}\label{sec3}

A computable approximation to the propagator can be defined by truncating its PCE approximation \eqref{PCE}. Because the PCE is a doubly infinite expansion in the directions of both the order $q$ of the Wick polynomials and the dimensionality $d$ of the truncated Gaussian probability space, one may naturally truncate the set of multi-indices as
\begin{eqnarray}
   \mathcal{I}_{d,q}:=\left\{\bm{\alpha}=(\alpha_{1},\alpha_{2},\cdots,\alpha_{d}),\, \alpha_{i}\in \mathbb{N}^{+}_{0}\,: \,|\bm{\alpha}|\le q \right\}.
\end{eqnarray}

The associated sparse truncation of the PCE approximation is then given as
\begin{eqnarray} 
   \textbf u(\bm x, t; \bm{\xi})=\sum_{\bm{\alpha}\in \mathcal {I}_{d,q}}\textbf{u}_{\bm{\alpha}}(\bm x, t) \mathcal{H}_{\alpha}(\xi).
\end{eqnarray}
Here $\bm{u}_{\bm{\alpha}}$ is the deterministic coefficients w.r.t the stochastic polynomial chaos.

Using the bi-orthogonal expansion and the graded lexicographic order of the multi-index $\bm{\alpha}$ in \eqref{solution_expan}, the solution $\mathbf{u}(\bm{x},t;\bm{\xi}): \Omega \times (0,T] \times D \to \mathbb{R}^2$ to our stochastic cylinder flow problem can be represented by the product of  stochastic and deterministic basis functions, i.e.

\begin{equation}
   \textbf u(\bm{x},t;\bm{\xi})=\sum_{i=1}^{p}\sum_{j=1}^{m} \bm{c}_{ij}(t)\psi_{j}(\bm{x})\mathcal{H}_{i}(\bm{\xi}),
   \label{solu-expansion}
\end{equation}
where for any fixed $t \in (0,T ]$, $\bm{c}_{ij}(t) \in \mathbb{R}^2$ are coefficients with respect to (w.r.t) the expansion, $\{\psi_{j}(\bm{x})\}_{j=1}^{m}$ are the deterministic basis for the physical domain and $\{\mathcal{H}_{i}(\bm{\xi})\}_{i=1}^{p}$ represents the  multivariate  Hermite polynomials that are orthonormal w.r.t. the joint probability measure $\rho(\bm{\xi})$ of a $d$-dimensional random variable $\bm{\xi}=(\xi_{1},\ldots,\xi_{d})$. Here $p=\frac{(q+d)!}{q!d!}$ is the cardinality of the set of multi-indices in \eqref{multi-indices}, where $\|\bm{\alpha} \|_{1} \leq q$,$\| \bm{\alpha}\|_{0} \leq d$, \cite{doostan2011non, gunzburger2014stochastic, schwab2011sparse, xiu2009fast}.

To determine the coefficients $\{\bm{c}_{ij}(t)\}_{1\le i\le p, 1\le j\le m}$ in expansion \eqref{solu-expansion}, the traditional CS method can be utilized when choosing the deterministic basis $\{\psi_{j}(\bm{x})\}_{j=1}^m$ as the finite element basis that are adopted to discretization w.r.t. $\bm{x} \in \Omega$. By enhancing the sparsity in \eqref{solu-expansion}, MDCS-PCE is a combination of CS and PCE through a limited number of sample solutions, to reduce the computational complexity further and provide a more efficient algorithm \cite{Liang}.

\subsection{Data-Driven compressive sensing approach}\label{sec3.1}

Now we introduce the data-driven compressive sensing method for solving SPDEs which the solutions are assumed to be sparse or compressible under the representation of some kind of polynomial basis.

Given a limited number of sample simulations, the data-driven compressive sensing method proposed in our previous work \cite{Liang}, constructs a problem-dependent basis for the expansion \eqref{solu-expansion}, which further enhances the sparsity in the representation of solution and hence improves the recovery accuracy. To be specific, we firstly construct the covariance function $\textnormal{Cov}_{\textbf{u}}(\bm{x},\bm{x}';t)$  of solution $\textbf{u}(\bm{x},t;\bm{\xi})$ using low-fidelity simulations which are less accurate but computationally cheaper than the high-fidelity ones  \cite{gunzburger2018}. Then, the widely studied Karhunen-Lo\`{e}ve analysis, based on the integral equation, is applied to extract the most dominant energetic modes as our data-driven basis functions.

Firstly, we recover the covariance function \eqref{Def-Cov} using the traditional compressive sensing method at low-fidelity mesh size when only a few observations are provided.
\begin{equation}
   \text{Cov}_{\textbf{u}}(\bm{x},\bm{x}';t)=\int_{D} \textbf u(\bm{x},t;\bm{\xi})\textbf{u}(\bm{x}',t;\bm{\xi})\rho(\bm{\xi}) d\bm{\xi}-\mathbb{E}_{\textbf{u}}(\bm{x},t) \mathbb{E}_{\textbf{u}}(\bm{x}',t).
   \label{Def-Cov}
\end{equation}
Here, $\mathbb{E}_{\textbf{u}}(\bm{x},t)$ denotes the expectation of $\textbf{u}(\bm{x},t;\bm{\xi})$. Given that solution $\textbf{u}(\bm{x},t;\bm{\xi})$ is sparse w.r.t. the basis functions in \eqref{solu-expansion}, the solution and hence its covariance can be accurately recovered using $n\ll p$ sample measurements via CS method \cite{doostan2011non, gwon2012compressive}.

By using the orthogonality of the multivariate Hermite polynomials $\{\mathcal{H}_{i}(\bm{\xi})\}_{i=1}^{p}$, we can derive the approximation of covariance function $\textnormal{Cov}_{\textbf{u}}(\bm{x},\bm{x}';t)$ as
\begin{equation}
   \textnormal{Cov}_{\hat{\textbf u}}(\bm{x},\bm{x}';t)=\sum_{i=2}^{p}\left(\sum_{j=1}^m\hat{\bm{c}}_{ij}(t)\psi_{j}(\bm{x})\right)\left(\sum_{j'=1}^m\hat{\bm{c}}_{ij'}(t)\psi_{j'}(\bm{x}')\right).
   \label{Cov-CS}
\end{equation}
where the coefficients $\hat{\bm{c}}_{ij}$ are determined by solving the basis pursuit (BP) problem \cite{Candes2006}.

Instead of using the standard finite element basis for the sparse representation of $\textbf u(\bm{x},t;\bm{\xi})$ in the spatial domain, we adopt a data-driven basis $\left\lbrace \bm{\phi}_{j}(\bm{x}) \right\rbrace_{j\geqslant 1}$ extracted from Karhunen-Lo\`{e}ve analysis \eqref{KLE-eqn-continuous} with a sequence of eigenpairs $\{\bm{\lambda}_{j},\bm{\phi}_{j}(\bm{x})\}_{j\geqslant 1}$, in which the kernel function is generated from \eqref{Cov-CS}.
\begin{equation}
   \int_{\Omega}\text{Cov}_{\textbf u}(\bm{x},\bm{x}')\bm{\phi}_{j}(\bm{x}')\,d\bm x'=\bm{\lambda}_{j}\bm{\phi}_{j}(\bm{x}).
   \label{KLE-eqn-continuous}
\end{equation}

Then the sparse solution of the problem \eqref{solu-expansion} can be derived by solving the convexified $\ell_1$-norm minimization problem \eqref{BP-Prob-DCS} associated with our data-driven basis.
\begin{equation}
   \tilde{\textbf{C}}=\operatorname*{arg\,min}_{\textbf{C}}\| \textbf{C}\|_1\ \ \textnormal{subject to}\ \ \| \bm{u}-\bm{\Psi}\textbf{C}\bm{\Phi} \|_2\leq \epsilon . 
   \label{BP-Prob-DCS}
\end{equation}
Here, for a fixed time $t$, $\bm{u}:=(\textbf{u}(\bm{x}_{i},\bm{\xi}_{j}))$ is the matrix of observations, $\textbf{C}:=(\tilde{\bm{c}}_{ij})$ is the coefficient matrix to be determined, $\bm{\Psi}:=(\mathcal{H}_{i}(\bm{\xi}_{j}))$ and $\bm{\Phi}:=(\bm{\phi}_{i}(\bm{x}_{j}))$ represent the spatial and stochastic information matrix respectively. 

The main MDCS-PCE algorithm is given as follows, we refer to \cite{Liang} for further details. 

\begin{algorithm}[htb]
   \caption{ $\langle$Multi-fidelity compressive sensing with data-driven basis$\rangle$}  
   \label{Algorithm-CS-KL-basis}
   \begin{algorithmic}
      \STATE {\textbf{Input}: coarse and fine mesh sizes ($h_{1}$ and $h_{2}$, respectively), stochastic collocation points. }
      \STATE{\textbf{Output}: compressed coefficients $\tilde{\textbf{C}}=(\tilde{\bm{c}}_{ij})$ w.r.t. data-driven basis for high-fidelity representation.}
      \STATE {1. Generate coarse mesh sample simulations and use the CS method to recover the covariance function \eqref{Cov-CS}.}
      \STATE {2. Solving Karhunen-Lo\`{e}ve decomposition problem \eqref{KLE-eqn-continuous} to compute data-driven basis $\{\bm{\phi}_j(\bm{x})\}_{j=1}^m$, then construct matrix $\bm{\Phi}=\left(\bm{\phi}_i(\bm{x}_j)\right)$ on fine mesh}.
      \STATE {3. Generate fine mesh sample simulations and solve the BP problem} \eqref{BP-Prob-DCS}\\
      \STATE {$$ 
              \tilde{\textbf{C}} =\operatorname*{arg\,min}_{\textbf{C}}\| \textbf{C}\|_1\ \ \textnormal{subject to}\ \ \| \bm{u}-\bm{\Psi}\textbf{C}\bm{\Phi} \|_2\leq \epsilon,$$}
      \STATE {\ \ \ \ then reconstruct the solution $\textbf u(\bm{x},t;\bm{\xi})$ with coefficients $\tilde{\textbf{C}}$, \textit{i.e.},}
      \vspace{2mm}
      \STATE {$$\tilde{\textbf u}(\bm{x},t;\bm{\xi}) = \sum_{i=1}^p\sum_{j=1}^m\tilde{\bm{c}}_{ij}(t) \phi_j(\bm{x}) \mathcal{H}_i(\bm{\xi}).$$ }
   \end{algorithmic}
\end{algorithm}

\begin{rmk}
   Given a stochastic flow field $\textbf{u}(\bm{x},t;\bm{\xi})$ satisfies \eqref{solu-expansion}, since the dominant solution modes do not vary substantially as time increases, we can divide the time domain $T$ into $N$ subintervals, instead of recovering covariance function at each time instance. Applying the same data-driven basis during each time interval may save the computational cost significantly.
\end{rmk}

\subsection{Simulation results using data-driven compressive sensing}\label{sec3.2}

In this section, we perform some numerical experiments to illustrate our numerical scheme and Algorithm \ref{Algorithm-CS-KL-basis}, with a focus on simulating deterministic and stochastic flows. In all the numerical simulations, the channel domain $\Omega=[0,1.2]\times[0, 0.5]$, the circular cylinder is centered at $(x_0, y_0)=(0.2, 0.25)$ and has diameter $d=0.1$, and the Reynolds number $Re=200$.

\subsubsection{ Deterministic cylinder flow }\label{sec3.2.1}

We first perform test on the $\theta-$scheme  presented in \cite{gunzburger2011optimal} on a deterministic problem by choosing $u_{in}(0,y,t;\bm{\xi})=\bar{u}_{in}=A y(H-y)$, where $A=70, H=0.5$, and $\phi(t,\bm{\xi})=\bar{\phi}=60\sin{(60\pi t)}$.

\begin{figure}[htb]
   \centering
   \includegraphics[width=3.3in,height=4.0in]{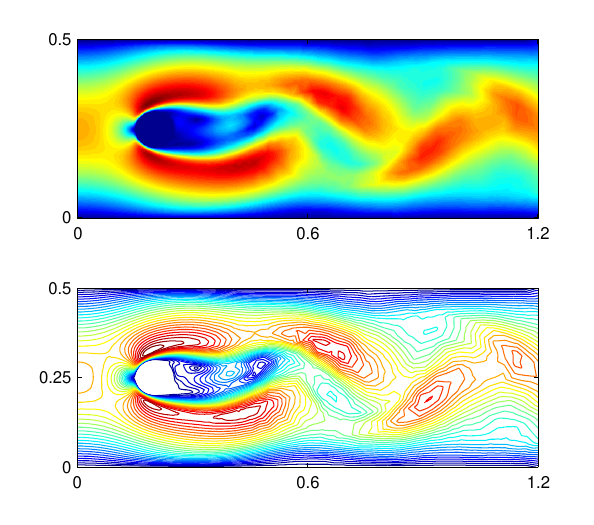}
   \caption{Scaled speed of the deterministic cylinder flow (top) and the associated contours (bottom) at $t=10$.\label{fig1}}
\end{figure}

\begin{figure}[htb]
   \centering
   \includegraphics[width=3.3in,height=1.8in]{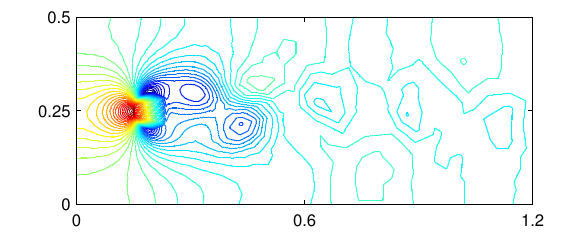}
   \caption{Pressure contours at $t=10$.\label{fig2}}
\end{figure}

Figure \ref{fig1} presents the instantaneous speed field $|\textbf u|:=\sqrt{u^{2}+v^{2}}$ of the deterministic cylinder flow and the associated contours at $t=10$; the pressure contours are given in Figure \ref{fig2}. The vortex street behind the circular cylinder is clearly visible.

\begin{figure}[htb]
   \centering
   \includegraphics[width=6.5cm]{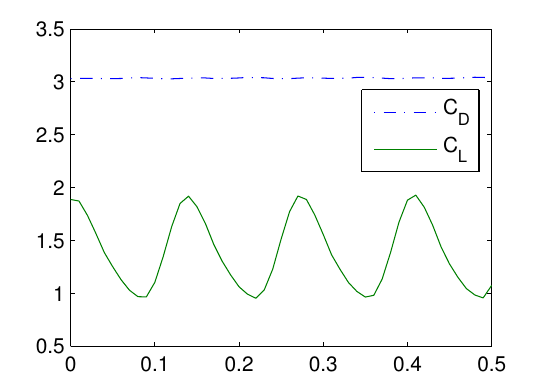}
   \caption{Time history of the drag coefficient $C_{D}$ and lift coefficient $C_L$ for the deterministic problem.
   \label{fig3}}
\end{figure}

It is well known that the dimensionless force on the cylinder is given by
\begin{gather}
   \textbf F=\big(F_{D}, F_{L} \big)=\int_{\Gamma_{c}}\left(p\textbf{n}-\frac{1}{Re}\frac{\partial\textbf{u}}{\partial \textbf{n}}\right)d\Gamma,
\end{gather}
where $F_D$ and $F_L$ are the drag and lift forces, respectively, and the corresponding dimensionless drag and lift coefficients are $C_{D}=2F_{D}$ and $C_{L}=2F_{L}$. Additionally, the \emph{Strouhal} number is defined as $S_{t} = fd/ \textbf{u}_{ave}$, where $f$ is the frequency of separation, i.e., the vortex shedding frequency, $d$ is the diameter of the cylinder, and here the average velocity $ \mathbf{u}_{ave}(t) = 2\mathbf{u}_{in}(0,H/2, t)/3$. We compare our coefficients $C_{D}$ and $C_{L}$ and the Strouhal number $S_{t}$ with the results obtained by the PARDISO solver of the COMSOL Multiphysics package. Table \ref{table2} indicates the agreement of the solutions obtained via the two different approaches. Time histories of $C_{D}$ and $C_{L}$ from $t=0$ to $0.5$ are plotted in Figure \ref{fig3}.

\begin{table}[htb]
    \label{table1}
    \centering
    \begin{tabular}{|c|| c| c| c|} 
    \hline 
    & $C_L$ & $C_D$ & $S_t$ \\
    [1ex]
   \hline  
   COMSOL Multiphysics package & 1.0898 & 2.9878 & 0.1146\\ 
   \hline
   semi-implicit scheme\cite{gunzburger2011optimal} & 1.0996 & 2.9890 & 0.1143\\
   \hline
   \end{tabular}
   \caption{For the deterministic simulation problem, comparison of the average lift coefficients $C_{L}$, the average drag coefficients $C_{D}$, and the {Strouhal} number $S_t$ obtained by COMSOL and the semi-implicit scheme \cite{gunzburger2011optimal} specialized to the deterministic case.
   \label{table2}}
\end{table}

\subsubsection{Stochastic cylinder flow }\label{sec3.2.2}

\begin{figure}[htb]
   \centering
   \includegraphics[width=3.5in,height=2.5in]{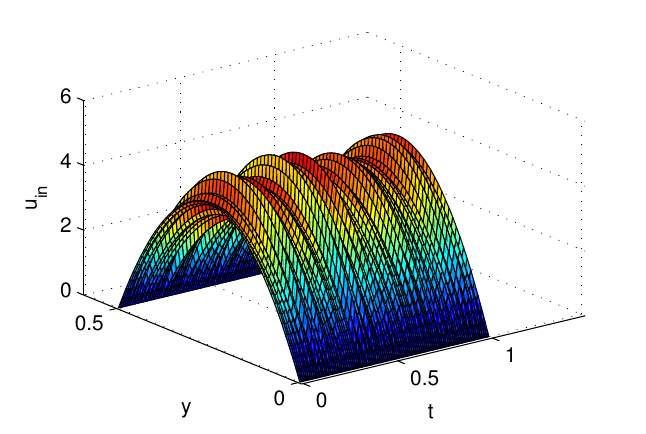}
   \caption{Time history of a realization of the stochastic inlet velocity $u_{in}(0,y,t;\bm{\xi})$ given in \eqref{stin} for $0\let\le 1$.}
   \label{fig4}
\end{figure}

\begin{figure}[htb]
   \centering
   \includegraphics[width=4.0in,height=1.5in]{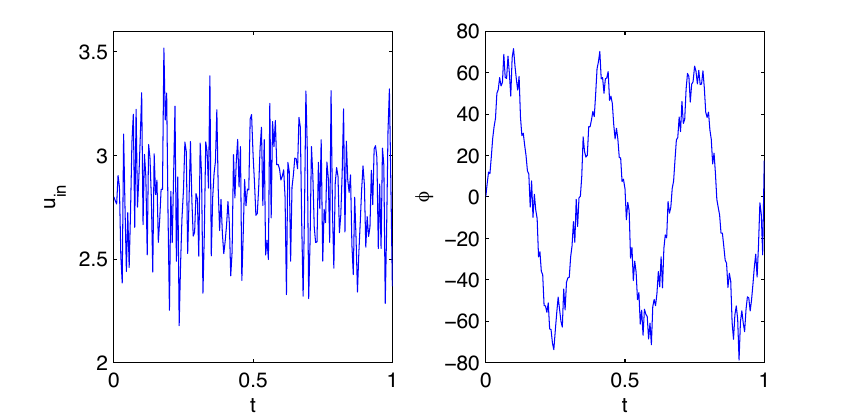}
   \caption{Time history of realizations of the stochastic inlet flow $u_{in}(0,y,t;\bm{\xi})$ at the point $(0, D/2)$ (left) and the stochastic angular velocity $\phi(t;\bm{\xi})$ (right) given in \eqref{stin}.}
   \label{fig5}
\end{figure}

We next consider applying our methodology to the simulation of the stochastic incompressible NS flows about a cylinder for which the random inputs $u_{in}(0,y,t;\bm{\xi})$ and $\phi(t;\bm{\xi})$ are both given. Specifically, in \eqref{bond1} and \eqref{gammac}, we respectively choose the stochastic inputs
\begin{equation}
   u_{in}(0,y,t;\bm{\xi})=(A+\sigma_{1}\dot\omega(t))y(H-y) \quad\mbox{and}\quad \phi(t;\bm{\xi})=\overline\phi(t)+\sigma_{2}\dot\omega(t),
   \label{stin}
\end{equation}
where $A=70$, $H=0.5$, $\sigma_{1}=2.0$, $\overline\phi(t)=60\sin(6\pi t)$, and $\sigma_{2}=1.0$. The \emph{Brownian} motion $\omega(t)$ is approximated by the direct truncation of \eqref{rep_Brownian}, i.e., we have
\begin{gather}
   \omega(t)\approx \omega_{K}(t):=\sum_{i=1}^{K}\xi_{i}\int_{0}^{t}\varphi_{i}(s)ds
   \label{app_noise}
\end{gather}
and the orthogonal basis of $L^2([0, T])$ is taken as
\begin{equation}
   \varphi_{i}(t)\d=\sqrt\frac{2}{{T}}\sin\left(\frac{(i-1/2)\pi t}{T}\right), i=1, 2, \ldots.
   \label{orth_basis}
\end{equation}
It is easily shown that
\begin{gather}
   \mathbb{E}\Big[\omega(t)-\sum_{i=1}^{K} \xi_{i}\int_{0}^{t}\varphi_{i}(s)ds\Big]^2\le \frac{t}{K},
\end{gather}
so that the error of the approximation \eqref{app_noise} is $\mathcal O (\sqrt{t})$. In our numerical test, we let $K=5$ and $\Delta t=10^{-3}$, then the errors incurred by truncating the expansions for the inputs at each time step from $t_{n}$ to $t_{n+1}=t_{n}+\Delta t$ are $\le 0.0063$. Figure \ref{fig4} and Figure \ref{fig5} show realizations of the inputs $u_{in}(0,y,t;\bm{\xi})$ and $\phi(t;\bm{\xi})$ for $0\le t \le 1$.

Due to the fact that the stochastic inputs only act on some portions of the boundary, a reasonable guess is that the second-order  PCE approximation would provide sufficient accuracy. We also note that the truncated error for the input noise approximation \eqref{app_noise} for $\Delta t=10^{-3}$ and $K=5$ is sufficiently small in our computation ($\le 0.0063$) so that to keep computational costs manageable, for the following numerical test of the data-driven compressive sensing approximation, we select $q=2, d=5$ as the upper bound for our sparse truncation. As a result, $p=21$ multivariate Hermite polynomials are used in the expansion \eqref{solu-expansion}, so the same number of polynomial chaos coefficients $\textbf{u}_{\bm{\alpha}}$ are generated.

\begin{table}
   \normalsize
   \centering
   \setlength{\tabcolsep}{3pt}
   \caption{\textnormal{Main eigenvalues and energy ratio for the Karhunen-Lo\`{e}ve decomposition of sparse scale $(h_{1})$ and fine scale $(h_{2})$ solutions at $t=0.5$}}\label{eigenvalue_t1}
   \begin{tabular}{c|c|c|c|c|c|c|c|c}
      \toprule
      $t=0.5$ & $\lambda_{1}$ & $\lambda_{2}$ & $\lambda_{3}$ & $\lambda_{4}$ & $\lambda_{5}$ & $\lambda_{6}$ & $\lambda_{7}$ & $\rho=\sum_{i=1}^{7}\lambda_{i} / \sum_{i=1}^{\infty} \lambda_{i}$\ \\
      \midrule
      $\text{Cov}_{u}^{h_{1}}$ & 4.81e-2 & 2.46e-2 & 2.09e-2 & 5.28e-3 & 1.06e-3 & 3.90e-4 & 2.41e-4 & 0.9982 \\
      \midrule
      $\text{Cov}_{u}^{h_{2}}$ & 6.06e-2 & 2.35e-2 & 1.38e-2 & 5.97e-3 & 1.83e-3 & 1.50e-3 & 5.16e-4 & 0.9948 \\
      \midrule
      $\text{Cov}_{v}^{h_{1}}$ & 4.97e-2 & 2.03e-2 & 1.40e-2 & 6.91e-3 & 5.27e-3 & 2.99e-3 & 9.12e-5 & 0.9988 \\
      \midrule
      $\text{Cov}_{v}^{h_{2}}$ & 9.48e-2 & 2.35e-2 & 3.47e-3 & 3.21e-3 & 1.92e-3 & 6.34e-4 & 4.47e-4 & 0.9968 \\
      \bottomrule
   \end{tabular}
\end{table}

\begin{table}
   \normalsize
   \centering
   \setlength{\tabcolsep}{3pt}
   \caption{\textnormal{Main eigenvalues and energy ratio for the Karhunen-Lo\`{e}ve decomposition of sparse scale $(h_{1})$ and fine scale $(h_{2})$ solutions at $t=1$}}
   \label{eigenvalue_t2}
   \begin{tabular}{c|c|c|c|c|c|c|c|c}
      \toprule
      $t=1$ & $\lambda_{1}$ & $\lambda_{2}$ & $\lambda_{3}$ & $\lambda_{4}$ & $\lambda_{5}$ & $\lambda_{6}$ & $\lambda_{7}$ & $\rho=\sum_{i=1}^{7}\lambda_{i} / \sum_{i=1}^{\infty} \lambda_{i}$\ \\
      \midrule
      $\text{Cov}_{u}^{h_{1}}$ & 9.86e-2 & 7.91e-2 & 4.32e-2 & 2.30e-3 & 1.17e-3& 8.77e-4 & 6.10e-4 & 0.9962 \\
      \midrule
      $\text{Cov}_{u}^{h_{2}}$ & 7.77e-2 & 5.18e-2 & 1.28e-2 & 2.51e-3 & 1.83e-3 & 8.53e-4 & 7.06e-4 & 0.9922 \\
      \midrule
      $\text{Cov}_{v}^{h_{1}}$ & 1.22e-1 & 9.60e-2 & 2.11e-3 & 1.26e-3 & 9.17e-4 & 5.77e-4 & 4.18e-4 & 0.9974 \\
      \midrule
      $\text{Cov}_{v}^{h_{2}}$ & 1.05e-1 & 2.01e-2 & 3.39e-3 & 2.16e-3 & 9.37e-4 & 5.32e-4 & 3.49e-4 & 0.9941 \\
      \bottomrule
   \end{tabular}
\end{table}
Firstly, to show the effectiveness of our MDCS-PCE algorithm, the distribution of eigenvalues of $\text{Cov}_{\mathbf{u}}(\bm{x},\bm{x}';t)$ for different mesh sizes is given in Table \ref{eigenvalue_t1} and Table \ref{eigenvalue_t2}, it can be seen that the main spectrum remains almost unchanged as the mesh gets refined, which ensures the feasibility of utilizing the low-cost coarse mesh sample simulations to generate the fine scale data-driven basis functions with little loss in accuracy. The data of energy ratio $\rho$ shows that the first 7 eigenvalues contain more than $99\%$ information of the solution $\mathbf{u}(\bm{x},t;\bm{\xi})$, which implies that the MDCS-PCE basis can compress the expansion coefficients extremely due to Theorem 3 in \cite{Liang}.

\begin{figure}[htb]
    \centering
    \includegraphics[width=3.2in,height=3.5in]{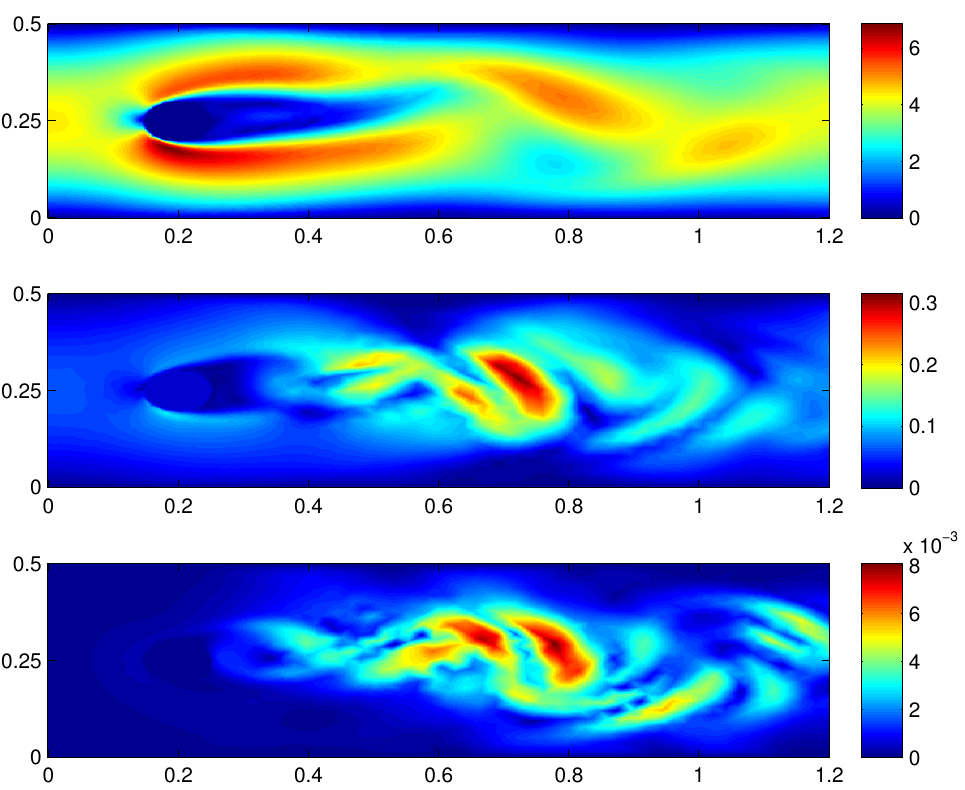}
   \caption{Scaled PCE coefficients $\mathbf{u}_{\bm{\alpha}}$ for the multi-indices $\bm{\alpha}=\textbf 0 $ (top) which corresponds to the mean value $\mathbb{E}[\textbf{u}]$; $\bm{\alpha}=(1,0,0,0,0)$ (middle); and $\bm{\alpha}=(0,0,0,0,2)$ (bottom) which is the last coefficient at $t=0.5$.}
   \label{fig_part3_image}
\end{figure}
\begin{figure}[htb]
    \centering
   \includegraphics[width=3.2in,height=3.5in]{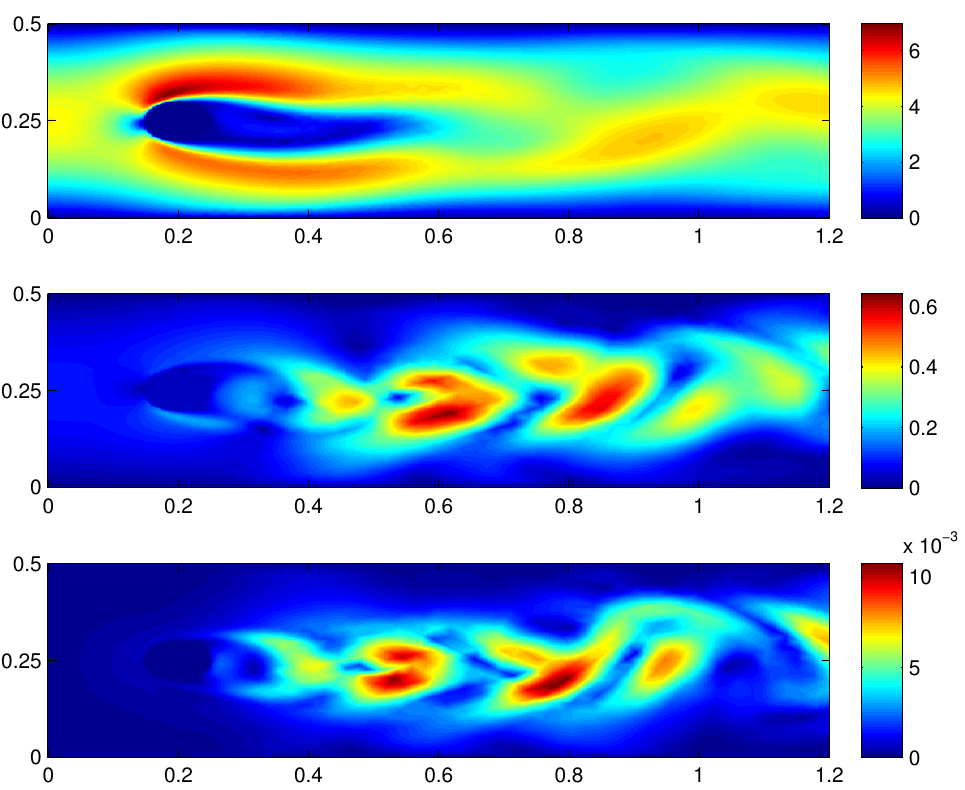}
   \caption{Scaled PCE coefficients $\textbf{u}_{\bm{\alpha}}$ for the multi-indices $\bm{\alpha}=\textbf 0 $ (top) which corresponds to the mean value $\mathbb{E}[\textbf{u}]$; $\bm{\alpha}=(1,0,0,0,0)$ (middle); and $\bm{\alpha}=(0,0,0,0,2)$ (bottom) which is the last coefficient at $t=1$.}
   \label{fig6}
\end{figure}
\begin{figure}[htb]
    \centering
    \includegraphics[width=3.2in,height=2.9in]{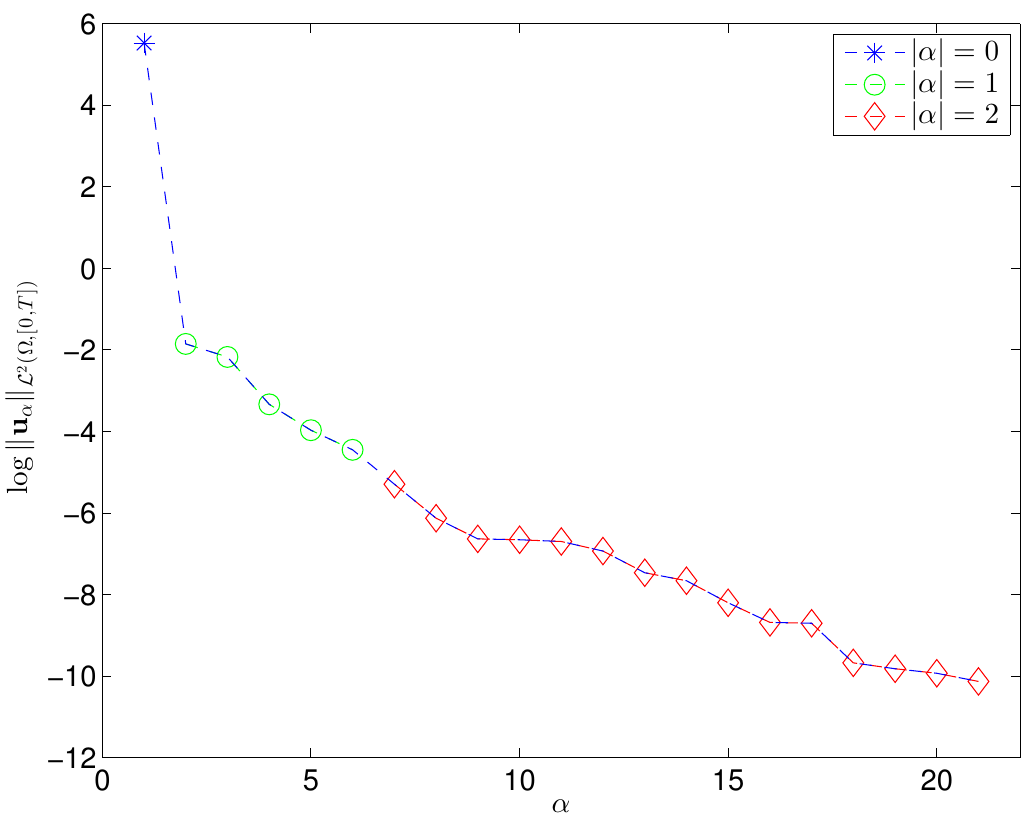}
   \caption{The time integral \eqref{l2n} of the $L^{2}(\Omega)$-norms of the 21 coefficients of the truncated PCE. The first 6 coefficients correspond to the Gaussian part of the solution (order $|\bm{\alpha}|\le 1$); coefficients 7 to 21 correspond to the second order Wick polynomials.}
   \label{fig7}
\end{figure}

Then numerical simulations of the MDCS-PCE method for three choices of the chaos coefficients are presented in Figure \ref{fig_part3_image} and \ref{fig6}. We also compute the $L^{2}(\Omega)$-norms of the coefficients $\textbf{u}_{\bm{\alpha}}$ over the time interval $[0,T]$ as given by
\begin{eqnarray}
   \int_{0}^{T}\|\textbf{u}_{\bm{\alpha}}\|^{2}_{L^{2}(\Omega)}dt,\qquad \bm{\alpha}\in \mathcal{I}_{5,2};
   \label{l2n} 
\end{eqnarray}
these are plotted in Figure \ref{fig7}. It can be seen that the PCE coefficients converge rapidly and the norms of the higher coefficients are $\le 10^{-4}$, this shows that the sparse truncation can provide sufficient accuracy for the simulation of an SNS flow, and also demonstrates the sparsity or compressible property of the solution, which is necessary for the application of our data-driven compressive sensing method.

\begin{figure}[htb]
    \centering
   	\includegraphics[width=3.3in,height=2.5in]{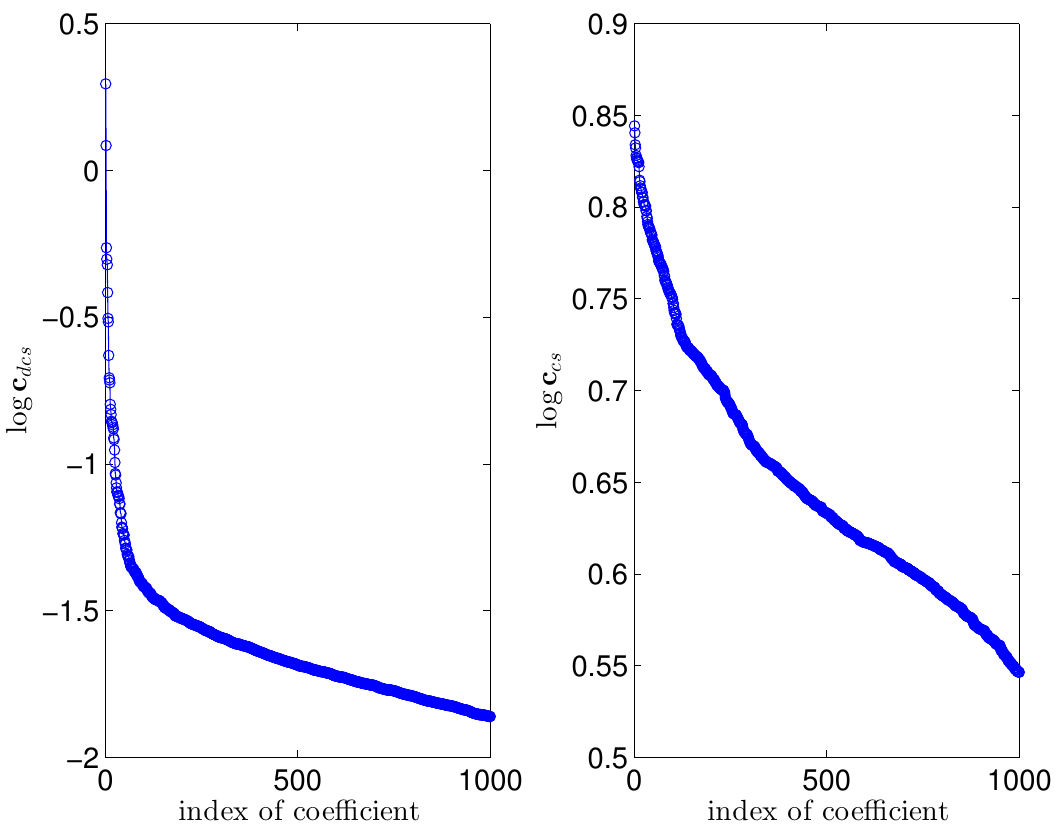}
   \caption{The first thousand big coefficients of both data-driven compressive sensing method solution and compressive sensing solution at $t=1$.}
   \label{fig9}
\end{figure}

Table \ref{eigenvalue_t1} and \ref{eigenvalue_t2} show the fast decay of the eigenvalues for the problem \eqref{KLE-eqn-continuous}, which demonstrates the feasibility of the MDCS-PCE method for further improving the sparse representation of the numerical approximation. The comparison between the traditional CS and our MDCS-PCE method are given in Table \ref{table-example-1-sparsity} and Figure \ref{fig9}. By comparing the number of larger coefficients in the expansion \eqref{solu-expansion}, it can be found that the MDCS-PCE method can significantly increase the sparsity of the coefficients than the traditional compressive sensing method.

\begin{table}
   \normalsize
   \centering
   \setlength{\tabcolsep}{3pt}
   \caption{\textnormal{Sparsity of the coefficient matrix $\bm{c} $ w.r.t. different methods and threshold values at $t=0.5$ and $t=1$. CS method: $\hat{\mathbf{C}}$, MDCS-PCE method: $\tilde{\mathbf{C}}$.}}
   \label{table-example-1-sparsity}
   \begin{tabular}{c|c|c||c|c|c}
      \toprule
      $t=0.5$ & $\#\{|\hat{\mathbf{C}}|\geq\tau\}$ & $\#\{|\tilde{\mathbf{C}}|\geq\tau\}$ & $t=1$ & $\#\{|\hat{\mathbf{C}}|\geq\tau\}$ & $\#\{|\tilde{\mathbf{C}}|\geq\tau\}$ \\
      \midrule
      $\tau=10^{-1}$ & 2814 & 21 & $\tau=10^{-1}$ & 3062 & 25 \\
      \midrule
      $\tau=10^{-2}$ & 13657 & 832 & $\tau=10^{-2}$ & 15694 & 1418 \\
      \midrule
      $\tau=10^{-3}$ & 26148 & 2313 & $\tau=10^{-3}$ & 28017 & 2383 \\
      \midrule
      $\tau=10^{-4}$ & 37514 & 2465 & $\tau=10^{-4}$ & 37492 & 2478 \\
      \bottomrule
      \end{tabular}
\end{table}

\section{Optimal boundary control of stochastic cylinder flow}\label{sec4}

In our stochastic flow control problem, the governing SNS equations \eqref{NS} with boundary conditions given by \eqref{bond1}--\eqref{donothingcond} are generally referred to as the \emph{state equations}. The control acts on the cylinder $\Gamma_{c}$ via adjusting the angular velocity $\phi$ in \eqref{gammac} so that we have a boundary control problem.

\subsection{Objective functionals}\label{sec4.1}  

There are several types of objective functionals for our boundary control problem, each serving its own set of control goals. Here we consider the \emph{Velocity tracking problem} as a typical example of stochastic flow control.

This problem reflects the desire to steer a candidate velocity field $\textbf{u}$ to a given target velocity field $\textbf{U}$, which can be deterministic or stochastic, over a certain period of time  by adjusting the fluid velocity along a portion of the boundary. Particularly, if the target velocity field is from a steady laminar flow, e.g., steady Stokes flow, then the unsteadiness of the controlled flow is expected to be minimized; if the target velocity field $\textbf {U}\equiv 0$, then the kinetic energy of the controlled flow is expected to be minimized; see, e.g., \cite{gunzburger2003perspectives, Hong2016The, graham1999optimal, Xu2017Flow, li2003optimal}.

The objective functional for the velocity tracking problem can naturally be chosen as
\begin{equation}
   \mathcal{K}(\textbf{u},\phi):=\mathbb{E}\left[\frac{1}{2}\int^{T}_{0}\Big((\textbf{u}-\textbf{U},\textbf{u}-\textbf{U})+ \beta|\phi|^{2}\Big)dt\right],
   \label{obj1}
\end{equation}
where the first term reflects the control goal. The second term provides a measurement of the control effort and expense, scaled by a weighting parameter $\beta\ge 0$, so that its inclusion in the functional has the purpose of ensuring that the control function lies in a reasonable range; see, e.g., \cite{tikhonov1977solutions, gunzburger2000adjoint, Kontoleontos2013Adjoint}. Often, the second term is also referred to as a \emph{penalty term} and $\beta$ as a \emph{penalty parameter}.

The boundary optimal control problem that we consider here is then simply stated as follows: find $\phi(t)$ such that the objective functional $\mathcal{K}(\textbf{u},\phi)$ is minimized subject to the SNS system \eqref{NS} and \eqref{bond1}--\eqref{donothingcond}.

\subsection{Optimization algorithm}\label{sec4.2}
We approximate the velocity $\textbf{u}$ by the MDCS-PCE method discussed in Algorithm \ref{Algorithm-CS-KL-basis}. This leads to approximations of the functional $\mathcal{K}$, we have the following approximations of $\mathcal{K}$ for the \emph{velocity tracking problem}.
\begin{equation}
   \begin{aligned}
      \mathcal{K}_{d,q}(\textbf{u},\phi) & \approx \frac{1}{2}\sum_{j=1}^{p}\int^{T}_{0}\left(\left(\textbf{u}_{j}-\textbf{U}_{j}, \textbf{u}_{j}-\textbf {U}_{j}\right) +\beta|\tilde{\phi}_{j}|^2\right)dt \\
      & =\frac{1}{2}\int^{T}_{0} \sum_{j=1}^{p}\left( \sum_{i=1}^{m}\left( \tilde{\bm{c}}_{ij}-\tilde{\bm{C}}_{ij} \right)^{2} + \beta |\tilde{\phi}_{j}|^{2} \right)dt.
   \end{aligned}
   \label{discfunc}
\end{equation}
To make notation uniform, we change the subscript `$\bm{\alpha}$'  to ` $j$ ', so here $\textbf{u}_{j}$ represents the expansion coefficients of the flow solution $\textbf{u}(\bm{x},t;\bm{\xi})$ while $\textbf{U}_{j}$ and $\tilde{\bm{C}}_{ij}$ are the coefficients of the given flow field $\textbf{U}(\bm{x},t;\omega)$, i.e. 
\begin{equation}
   \textbf{U}(\bm{x},t;\bm{\xi})=\sum_{j=1}^{p}\textbf{U}_{j} \mathcal{H}_{j}(\bm{\xi}), \quad \quad \textbf{U}_{j}=\sum_{i=1}^{m} \tilde{\bm{C}}_{ij}(t)\bm{\phi}_{i}(\bm{x}),
   \label{U_expansion}
\end{equation}
and $\tilde{\phi}_{j}(t)$ is the expansion coefficient for the control function $\phi$, i.e.,
$$
\phi(t;\bm{\xi})=\sum_{j=1}^{p}\tilde{\phi}_{j}\mathcal{H}_{j}(\bm{\xi}).
$$
Our approximate optimal control problem based on MDCS-PCE is then to determine the coefficients $\{\tilde{\phi}_{j}(t)\}_{j=1}^{p}$ such that $\d\mathcal{K}$ is minimized subject to the equation \eqref{NS} being satisfied.
%
%
We introduce the Lagrange multipliers (or adjoint) functions,
$$
\pmb\theta=(\theta^{(1)}, \theta^{(2)}) \quad  \theta^{(i)}(\bm{x},t;\bm{\xi}) \in \mathcal{L}^{2}(D;0,T;H^{1}(\Omega)), \quad i=1,2.
$$ 

Then we define the Lagrangian functional

\begin{equation}
   \begin{array}{l}
      \mathcal{M}(\textbf{u};\phi;\pmb\theta):=\mathcal{K}\big(\textbf{u};\phi\big) - \displaystyle\mathbb{E}\int^{T}_{0} (\pmb\theta, \frac{\partial \textbf{u}}{\partial t}-\frac{1}{Re}\Delta\textbf{u} + (\textbf{u}\cdot \nabla)\textbf{u}+ \nabla p \big)dt.
      \end{array}
\end{equation}

As what we need in our MDCS-PCE method is some sample simulation results for both the \emph{state} and the \emph{adjoint } equations, so for simplicity, we only need to consider the deterministic form. Setting the first variation of $\mathcal{M}$ w.r.t. $\pmb\theta$ and $\textbf{u}$ to zero respectively results in the weak form of the state system \eqref{weak_form} and the weak form of the \emph{adjoint} equations \eqref{adjeqeq},
\begin{equation}
   \begin{aligned}
      \left( \frac{\partial\pmb\theta}{\partial t}, \textbf{v}\right)-a(\pmb\theta, \textbf{v})-c\left(\textbf{v}, \textbf{u},\pmb\theta\right)+c\left(\textbf{u}, \textbf{v},\pmb\theta\right)&=-\left(\mathbf{u}-\mathbf{U}, \textbf{v}\right)\qquad \forall \textbf{v} \in H^{1}(\Omega) \\
      \pmb\theta\big|_{t=T} & = \textbf{0}.
   \end{aligned}
   \label{adjeqeq}
\end{equation}
Note that \eqref{adjeqeq} is linear in adjoint variable $\bm{\theta}$.

Consequently, the variation of the Lagrangian functional takes the form
\begin{eqnarray}
   \delta \mathcal{M} = \mathbb{E}\Big(\d\int_{0}^{T}\left(\phi-\frac{1}{Re}\left(-(y-y_{0})\frac{\partial \theta^{(1)}}{\partial\textbf{n}}+(x-x_{0})\frac{\partial\theta^{(2)}}{\partial\textbf{n}}\right),\delta\phi\right)_{\Gamma_{c}}dt\Big),
\end{eqnarray}
where $(\cdot, \cdot)_{\Gamma_{c}}$ is the $L^{2}$ inner product over boundary $\Gamma_{c}$.

If we let
\begin{eqnarray}
   \delta\phi=-\lambda\left(\phi-\frac{1}{Re}\left(-(y-y_{0})\frac{\partial \theta^{(1)}}{\partial\textbf{n}}+(x-x_{0})\frac{\partial\theta^{(2)}}{\partial\textbf{n}} \right)\right)\quad \lambda>0,
   \label{decent_dir}
\end{eqnarray}
on $\Gamma_{c}$, then $\delta\mathcal M<0$, which ensures, for sufficiently small step size $\lambda>0$, that the values of the objective functional $\mathcal {M}$ are reduced.

Based on the above discussion, we define the following gradient-based optimization algorithm.
\begin{algorithm}
   \caption{$ \langle$ Optimization algorithm based on the MDCS-PCE method$\rangle$}
   \label{Algorithm-Optimal-control}
   \begin{algorithmic}
      \STATE{1. Initialize the optimization.}
      \STATE{The optimization loop is initialized as follows.
      \begin{description}
         \item{a)} Depending upon the accuracy desired for the PCE expansion and the sparsity property, select an upper bound $q$ and the stochastic dimension $d$ to truncate the set of multi-indices $\mathcal{I}$, which we shall denote by $\mathcal I_{d,q}$.
         \item{b)} Select a positive step size parameter $\lambda$ and a tolerance parameter $\varepsilon$.
         \item{c)} Select finite element spaces $V^{h}$ and $Q^{h}$ for the velocity and pressure approximations and select a mesh size $h$. Select a time step $\Delta t$ .
         \item{d)}  Select a deterministic initial condition $\textbf{u}_{0}$ and stochastic boundary conditions $$ u_{in,d}(0,y,t;\bm{\xi})= \sum_{i\in\mathcal {I}_{d,q}}u_{in,i}(0,y,t)\mathcal{H}_{i}(\xi)$$ and $$\phi_{d}(t;\bm{\xi})=\sum_{i\in\mathcal{I}_{d,q}}\phi_{i}(t) \mathcal{H}_{i}(\xi).$$ 
         \item{e)} Using the data-driven compressive sensing method determine the coefficients $\{\textbf{u}_{i}\}_{i\in\mathcal{I}_{d,q}}$ from Algorithm \ref{Algorithm-CS-KL-basis}.
         \item{f)} Determine the value ${\mathcal{K}}_{d,q}$ of the functional \eqref{discfunc}.
    \end{description}}
    \STATE {2. Optimization loop.}
    \STATE {For $\ell=1,\ldots$,
    \begin{description}
       \item{a)} From $\{\phi_{d}^{\ell-1}\}$ and $\{\textbf{u}_{d,q}^{\ell-1}\}$, determine the set of solutions $\{\pmb\theta_{i}^{\ell-1}\}_{i \in \mathcal{I}_{d,q}}$ of the PCE-based adjoint system \eqref{adjeqeq}.
       \item{b)} From $\{\lambda^{\ell-1}\}$, $\{\phi_{d}^{\ell-1}\}$, and $\{\pmb\theta_{i}^{\ell-1}\}_{i\in\mathcal{I}_{d,q}}$, determine the set of steps $\{\delta\phi_{d}^{\ell}\}$ from \eqref{decent_dir}.
       \item{c)} From $\{\phi_{d}^{\ell-1}\}$ and $\{\delta\phi_{d}^{\ell}\}$, update the values of the chaos coefficients of boundary data from
       $$
       \phi^{\ell}_{d} = \phi^{\ell-1}_{d} + \delta\phi^{\ell-1}_{d}.
       $$
       \item{d)} From $\{\phi_{d}^{\ell}\}$, determine the set of coefficients $\{\textbf u_i^{\ell}\}_{i\in\mathcal I_{d,q}}$ using Algorithm \ref{Algorithm-CS-KL-basis}.
       \item{e)} From $\{\phi_{d}^{\ell}\}$ and $\{\textbf u_{i}^{\ell}\}_{i\in\mathcal{I}_{d,q}}$, determine the value of the functional ${\mathcal{K}}_{d,q}^{\ell}$ from \eqref{discfunc}.
       \item{f)} If ${\mathcal{K}}_{d,q}^{\ell}<{\mathcal{K}}_{d,q}^{\ell-1}$ and $({\mathcal{K}}_{d,q}^{\ell-1}-{\mathcal{K}}_{d,q}^{\ell})/{\mathcal{K}}_{d,q}^{\ell-1}\le\varepsilon$, stop and accept the result.
       
       If ${\mathcal{K}}_{d,q}^{\ell}<{\mathcal{K}}_{d,q}^{\ell-1}$ and $({\mathcal{K}}_{d,q}^{\ell-1}-{\mathcal{K}}_{d,q}^{\ell})/{\mathcal{K}}_{d,q}^{\ell-1}>\varepsilon$, go to Step 2a, incrementing the iteration counter $\ell$.
       
       If ${\mathcal{K}}_{d,q}^{\ell}>{\mathcal{K}}_{d,q}^{\ell-1}$, set $\lambda^{\ell-1}=\lambda^{\ell-1}/2$ and go to Step 2b.
       \end{description}}
   \end{algorithmic}
\end{algorithm}

We remark here that other gradient-based methods such as conjugate gradient or quasi-Newton methods are possible; see, e.g., \cite{gunzburger2000adjoint}.

\subsection{Computational experiments}\label{sec4.3}

To validate our optimization algorithm, numerical experiments with the same settings in Section \ref{sec3.2.2} are performed, i.e., $\Omega=[1.2]\times[0.5]$, the diameter of the cylinder $d=0.1$, the center $(x_{0}, y_{0})=(0.2, 0.25)$, and the Reynolds number $Re=200$, and the stochastic boundary data $u_{in}(0,y,t; \bm{\xi})=(70+2\dot{\omega}(t))(H-y)y$. The truncated tracking functional $\mathcal{K}_{d,q}\big(\textbf{u};\phi \big)$ is given by \eqref{discfunc}.

The initial condition $\textbf{u}_{0}$ is deterministic and chosen from the solution of the example in Section \ref{sec3.2.1} at $t=0.5$. In our numerical test, to validate the suppression of vortex shedding by tracking a steady velocity field, we select the target velocity $\textbf{U}$ from a deterministic NS flow past a stationary cylinder with no vortex shedding occurring. We let the initial guess for the control function $\bar{\phi}(t)=60 \sin(6\pi t)$, then the control is realized by adjusting $\mathbb{E}[\phi]$ using the optimization algorithm above. In fact, from Figure \ref{fig7} and Figure \ref{fig9}, we can observe that the stochastic part of the PCE coefficients ($|\mathbf{u}_{\bm{\alpha}}|, |\bm{\alpha}| >1$) are relatively small, which implies the applicability of just adjusting the expectation of $\phi$ for our stochastic optimization problem. 

\begin{figure}[htb]
   \centering
   \includegraphics[width=3.2in,height=3.7in]{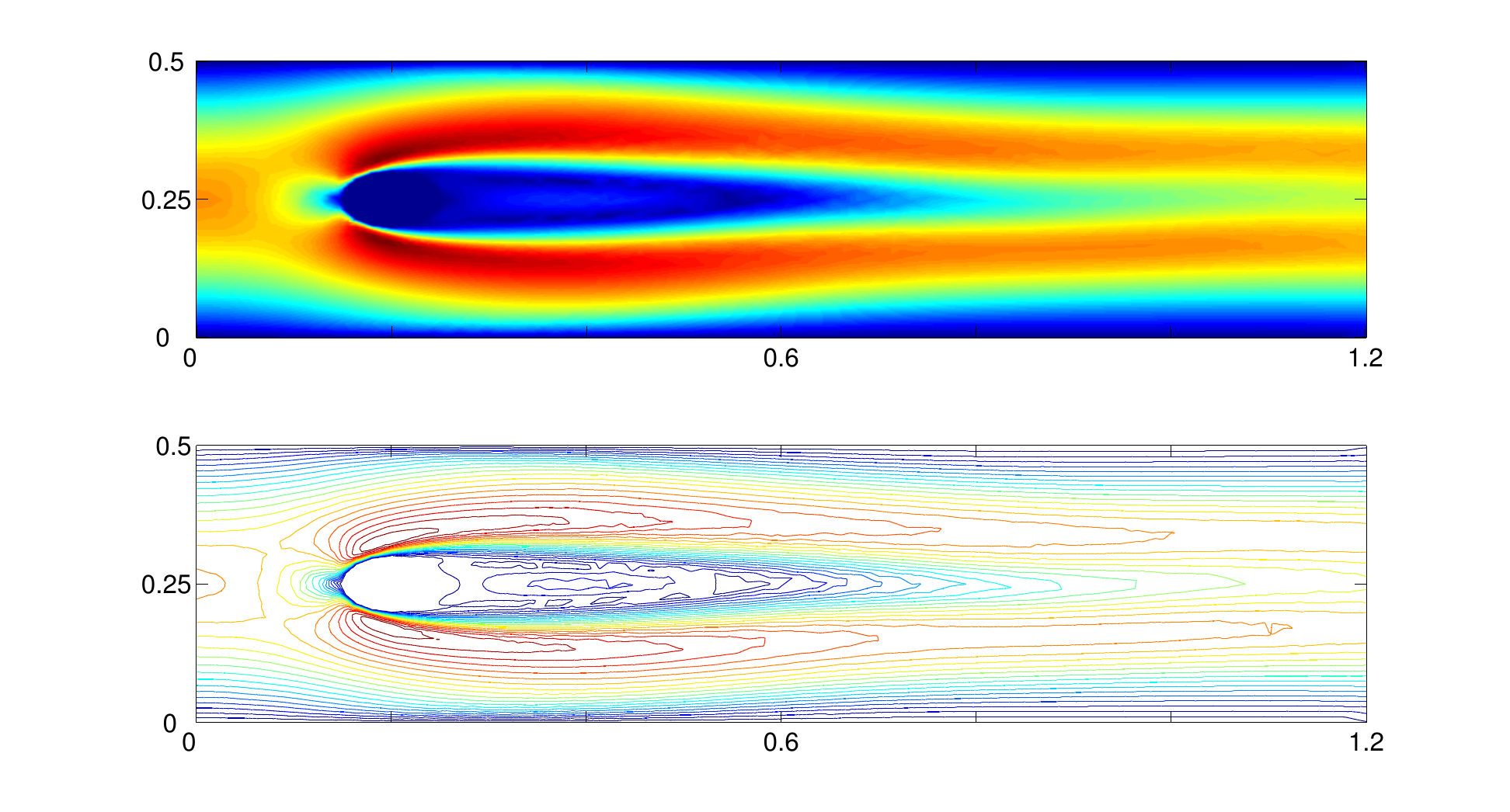}
   \caption{ Scaled average velocity of the deterministic flow past a stationary cylinder (top) and the associated contours (bottom) at $t=0.5$.
   \label{fig_U_tar}}
\end{figure}

\begin{figure}[htb]
    \centering
    \includegraphics[width=4.0in]{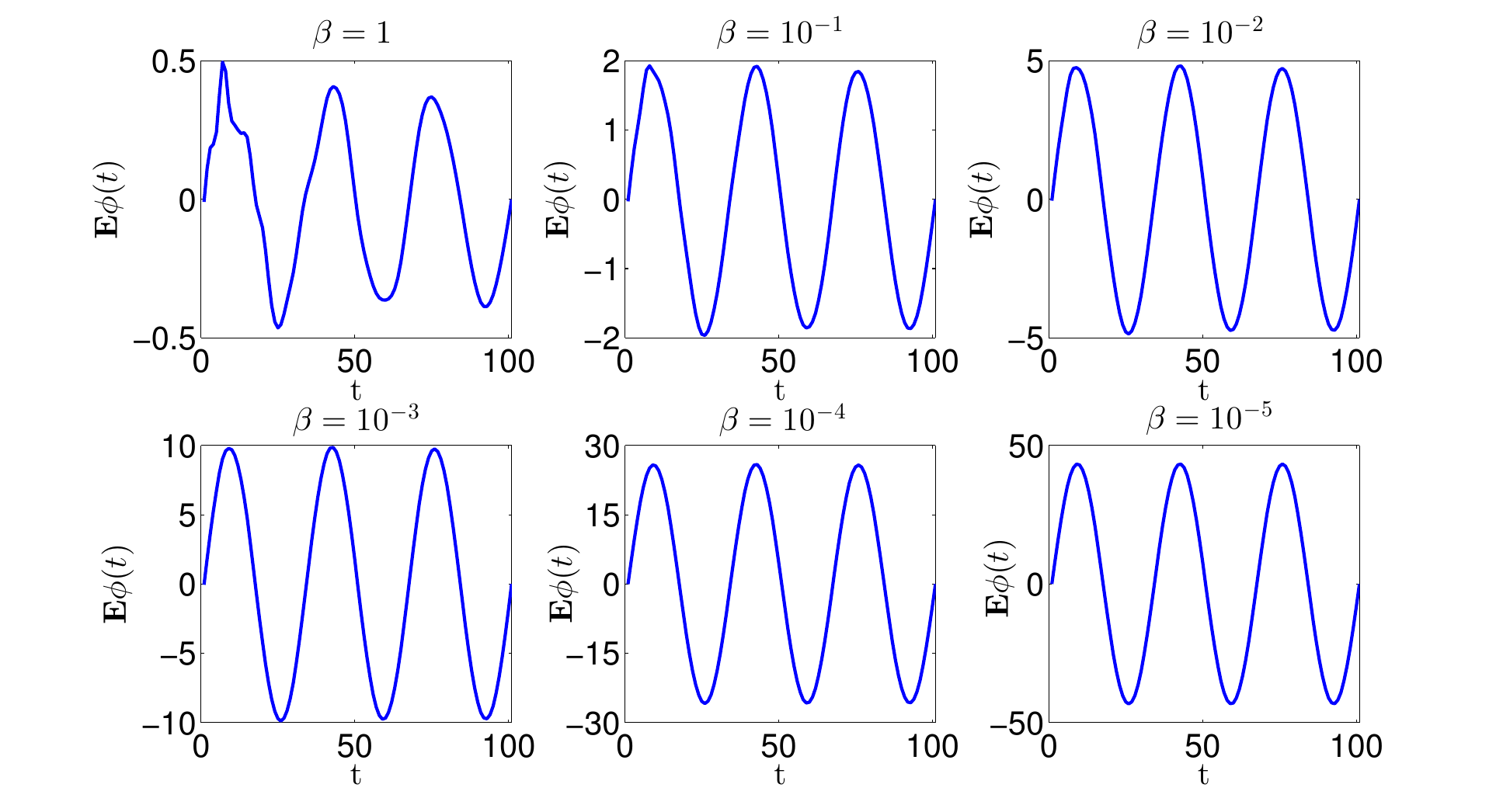}
   \caption{Optimized control actions for different $\beta$. Note that smaller $\beta$ corresponds to a larger range of optimal control action.}
   \label{fig_phi}
\end{figure}

\begin{figure}[htb]
    \centering
    \includegraphics[width=5.2in]{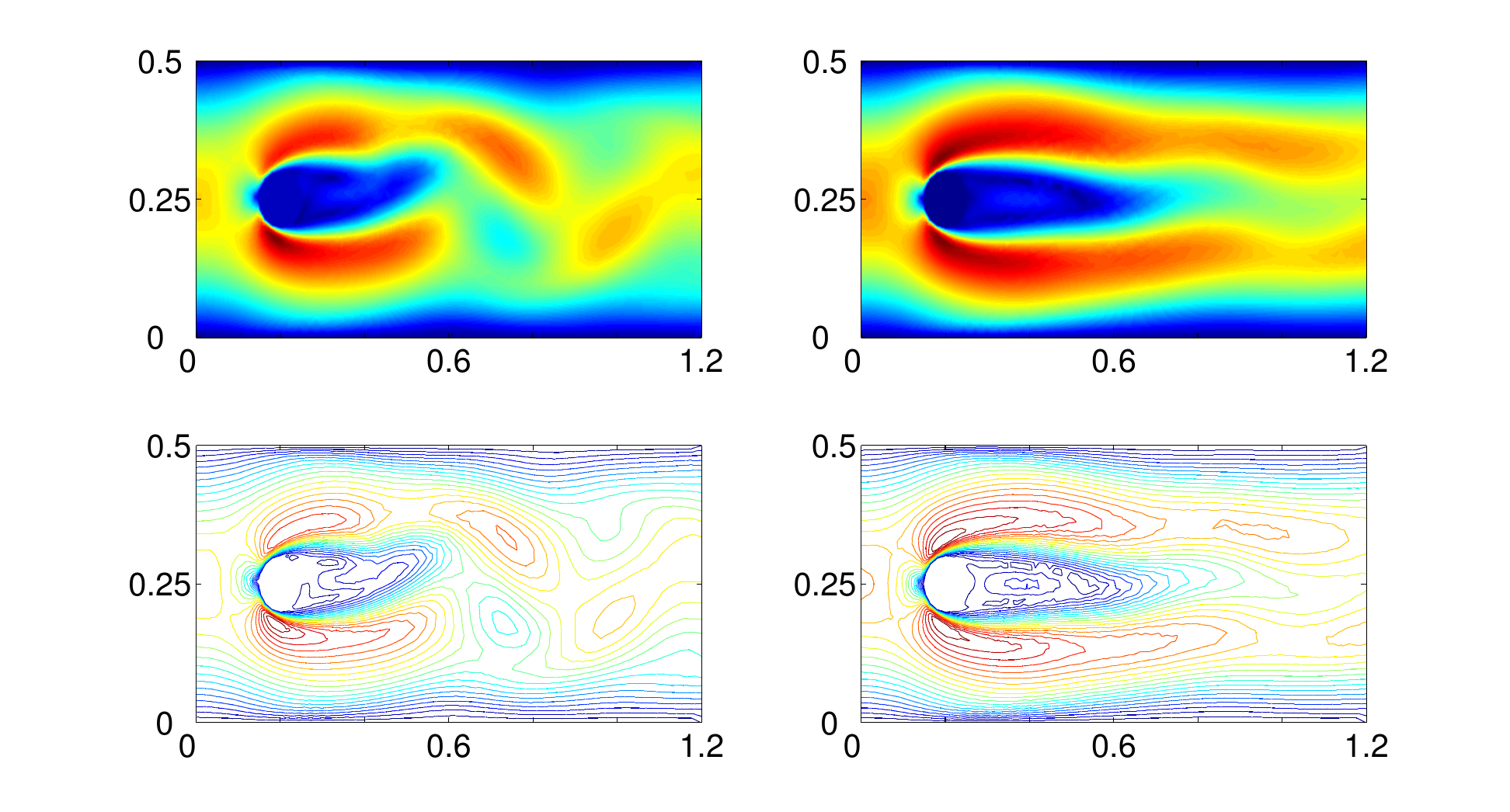}
   \caption{Scaled average velocity field and the associated contours of the optimized stochastic cylinder flow as $\beta=10^{-3}$ (first column), and $\beta=10^{-5}$ (second column) at $t=1.$}
   \label{fig_opt}
\end{figure}

Figure \ref{fig_U_tar} presents the scaled average velocity field of the deterministic cylinder flow past a stationary cylinder at $t=0.5$ when the vortex shedding is not established, we choose this state as the target velocity state $\textbf{U}$. Table \ref{table_opt} shows the control effect by assigning different values of parameter $\beta$, the associated control actions are presented in Figure \ref{fig_phi}, which certificate our assertion regarding the penalty term in Section \ref{sec4.1}: small $\beta$ corresponds to larger range of $\phi$, or larger weight of the control represented by the first term in $\mathcal{K}(\mathbf{u},\phi)$. According to Figure \ref{fig_opt}, as $\beta=10^{-5}$, the controlled cylinder flow attains its steady state, i.e., no vortex shedding occurs; however, if $\beta \geq 10^{-3}$, the range of control $\phi$ is not large enough to fully stabilize the flow. 

\begin{table}[htb]
    \centering
   \begin{tabular}{|c||c|c|c|c|}
   \hline 
   $\beta$  & 1 & $10^{-1} $ & $10^{-3}$ & $10^{-5}$ \\
   \hline 
   ini.$\mathcal{K}_{d,q}$ & $900.89$ & 90.389 & $1.2727$ & 0.5324\\
   \hline 
   opt.$\mathcal{K}_{d,q}$  &  0.4271 & 0.4094  & $0.3520$ & 0.0639 \\
   \hline
   \end{tabular}
   \caption{The initial and optimized values of objective functional \eqref{discfunc} for the boundary control of stochastic cylinder flow for different values of the penalty parameter $\beta$.} 
   \label{table_opt}
\end{table}

\begin{table}[htb]
    \centering
   \begin{tabular}{|c||c|c|c|}
   \hline 
   $\beta=10^{-5}$ & MDCS-PCE & CS & MC \\
   \hline 
   number of simulations & 80 & 120 & 2000 \\
   \hline 
   ini. $\mathcal{K}_{d,q}$ & 0.5324 & 0.5492 & 0.6037 \\
   \hline 
   opt. $\mathcal{K}_{d,q}$ & 0.0639 & 0.0721 & 0.0644  \\
   \hline
   \end{tabular}
   \caption{The comparison of the number of simulations that are needed to achieve at same optimal control goals for the three methods.}
   \label{table_compare}
\end{table}
In Table \ref{table_compare}, the number of simulations or measurements needed to achieve the same control goals for the three methods are stated, as it can be seen, the MDCS-PCE method than the traditional CS method, and both of the two methods are much more effective than the MC method.

\section{Conclusion}\label{sec5}           
        
In this paper, we apply the MDCS-PCE method to optimal control problems for stochastic cylinder flows. By establishing a problem-dependent basis using the Karhunen-Lo\`{e}ve analysis and providing a more sparsified expression for the numerical sample solutions, our multi-fidelity MDCS-PCE method can reduce the computational costs drastically by reducing the number of measurements without losing accuracy, and thus providing an efficient and accurate algorithm for stochastic flow control problems. 

For problems in which random inputs act in a more widespread manner, the MDCS-PCE method can also work well since the sample solutions always have sparsity as the degrees of freedom increase. Indeed, our algorithm can be used directly for other types of stochastic control problems. In our future work, we would like to combine the MDCS-PCE method with the reduced-order modeling method for solving stochastic control problems.


\bibliographystyle{ieeetr}


\end{document}